\documentclass{gtart_h}
\def\ifplaintex{\expandafter\ifx\csname documentclass\endcsname\relax}

\def\gtp{{\mathsurround=0pt\it $\cal G\mskip-2mu$eometry \&\ 
$\cal T\!\!$opology $\cal P\!$ublications}}  % GT publications

\def\recd{{\small Received:\qua\receiveddate\ifx\reviseddate\relax
\else\qquad Revised:\qua\reviseddate\fi\par}} 

\def\lognumber#1{\def\thelognumber{#1}}
\def\volumenumber#1{\def\thevolumenumber{#1}}
\def\volumeyear#1{\def\thevolumeyear{#1}}
\def\papernumber#1{\def\thepapernumber{#1}}
\def\pagenumbers#1#2{\def\startpage{#1}\def\finishpage{#2}}
\def\published#1{\def\publishdate{#1}}

\def\received#1{\def\receiveddate{#1}}

\def\accepted#1{\def\accepteddate{#1}}
\def\asciititle#1{\def\theasciititle{#1}}
\def\covertitle#1{\def\thecovertitle{#1}}

\def\asciiaddress#1{\def\theasciiaddress{#1}}

\long\def\asciiabstract#1{\long\def\theasciiabstract{#1}}

\let\\\par\let\thelognumber\relax\let\thevolumenumber\relax
\let\thepapernumber\relax\let\thevolumeyear\relax\let\startpage\relax
\let\finishpage\relax\let\publishdate\relax\let\receiveddate\relax
\let\reviseddate\relax\let\accepteddate\relax\let\theasciititle\relax
\let\thecovertitle\relax\let\theasciiauthors\relax\let\theasciiaddress\relax
\let\theasciiabstract\relax

\let\theasciiemail\relax

\ifplaintex
\font\logobig=cmssbx10 scaled 3836
\font\logomed=cmssbx10 scaled 2557
\else
\font\logobig=cmssbx10 scaled 4200
\font\logomed=cmssbx10 scaled 2800
\fi

\long\def\makeagttitle{   %%% start of definition of \makeagttitle
\count0=\startpage
\agt\hfill      %   Journal title (top left) 
\hbox to 45truept{\vbox to 0pt{\vglue -13truept{\logomed A\kern -.37em{\logobig 
T}\kern -.38em G}\vss}\hss}
\break
{\small Volume \thevolumenumber\ (\thevolumeyear)
\startpage--\finishpage\nl
Published: \publishdate}

\vglue .25truein

{\parskip=0pt\leftskip 0pt plus
1fil\def\\{\par\smallskip}{\Large\bf\thetitle}\par\medskip} \vglue
0.05truein

{\parskip=0pt\leftskip 0pt plus 1fil\def\\{\par}{\sc\theauthors}
\par\medskip}%
 
\vglue 0.03truein 

{\small\leftskip 25truept\rightskip 25truept{\bf Abstract}\stdspace\theabstract

{\bf AMS Classification}\stdspace\theprimaryclass
\ifx\thesecondaryclass\relax\else; \thesecondaryclass\fi\par
{\bf Keywords}\stdspace \thekeywords\par}\vglue 7truept

}   %%%% end of definition of \makeagttitle

\ifplaintex
\font\phead=cmsl9 scaled 950
\font\pnum=cmbx10 scaled 913
\font\pfoot=cmsl9 scaled 950
\headline{\vbox to 0pt{\vskip -4.5mm\line{\small\phead\ifnum
\count0=\startpage ISSN 1472-2739 (on-line) 1472-2747 (printed)
\hfill {\pnum\folio}\else\ifodd\count0\def\\{ }% 
\ifx\theshorttitle\relax\thetitle\else\theshorttitle\fi\hfill{\pnum\folio}
\else\def\\{ and }{\pnum\folio}\hfill\ifx\theshortauthors\relax\theauthors
\else\theshortauthors\fi\fi\fi}\vss}}
\footline{\vbox to 0pt{\vglue 0mm\line{\small\pfoot\ifnum\count0=\startpage
\copyright\ \gtp\hfill\else
\agt, Volume \thevolumenumber\ (\thevolumeyear)\hfill\fi}\vss}}
\else
\font=cmsl9 scaled 1050
\font\lnum=cmbx10 
\font=cmsl9 scaled 1050
\makeatletter
\def\@oddhead{{\small\lhead\ifnum\count0=\startpage ISSN 1472-2739 
(on-line) 1472-2747 (printed)\hfill {\lnum\number\count0}\else\ifodd\count0
\def\\{ }\ifx\theshorttitle\relax \thetitle \else\theshorttitle\fi\hfill
{\lnum\number\count0}\else\def\\{ and }{\lnum\number\count0}
\hfill\ifx\theshortauthors\relax 
\theauthors\else\theshortauthors\fi\fi\fi}}\def\@evenhead{\@oddhead}
\def\@oddfoot{\small\lfoot\ifnum\count0=\startpage\copyright\ \gtp\hfill\else
\agt, Volume \thevolumenumber\ (\thevolumeyear)\hfill\fi}
\def\@evenfoot{\@oddfoot}
\makeatother
\fi
\let\maketitlepage\makeagttitle
\let\makeshorttitle\maketitlepage
\let\maketitle\maketitlepage

\newwrite\gtoutfile
\long\gdef\makeheadfile{  %%% start of definition of \makeheadfile
{\def\\{, }\def\s{ }
\immediate\openout\gtoutfile head.xxx
\immediate\write\gtoutfile{Proxy-for: \ifx\theasciiauthors\relax
\theauthors\else\theasciiauthors\fi\s<\ifx\theasciiemail\relax\theemail\else\theasciiemail\fi>}
\immediate\write\gtoutfile{\noexpand\\}
\immediate\write\gtoutfile{Authors: \ifx\theasciiauthors\relax
\theauthors\else\theasciiauthors\fi}
{\def\\{ }\immediate\write\gtoutfile{Title: \ifx\theasciititle\relax
\thetitle\else\theasciititle\fi}}
\immediate\write\gtoutfile{Subj-class: GT or SG, GR etc}
\immediate\write\gtoutfile{MSC-class: \theprimaryclass\ifx\thesecondaryclass\relax\else, \thesecondaryclass\fi}
\immediate\write\gtoutfile{Journal-ref: Algebr. Geom. Topol. \thevolumenumber\s
(\thevolumeyear) \startpage-\finishpage}
\immediate\write\gtoutfile{Comments: Published by Algebraic and
Geometric Topology at}
\immediate\write\gtoutfile{\s\s\s  http://www.maths.warwick.ac.uk/agt/AGTVol\thevolumenumber/agt-\thevolumenumber-\thepapernumber.abs.html}
\immediate\write\gtoutfile{\noexpand\\}
\immediate\write\gtoutfile{}
\ifx\theasciiabstract\relax
\immediate\write\gtoutfile{\theabstract}\else
\immediate\write\gtoutfile{\theasciiabstract}\fi
\immediate\write\gtoutfile{}
\immediate\write\gtoutfile{\noexpand\\}
\immediate\write\gtoutfile{}
\immediate\closeout\gtoutfile}}  %%% end of definition of \makeheadfile

\def\maketitlepage{\makeagttitle\makeheadfile}
\let\makeshorttitle\maketitlepage
\let\maketitle\maketitlepage

\lognumber{19}
\volumenumber{5}
\volumeyear{2005}
\papernumber{19}
\pagenumbers{419}{442}
\received{4 December 2003} 
\accepted{6 May 2005}
\published{24 May 2005}

\usepackage{amssymb,amsmath}
\usepackage[all]{xy}

\def \R{\mathbb{R}}
\def \N{\mathbb{N}}
\def \Z{\mathbb{Z}}
\def \Q{\mathbb{Q}}
\def \C{\mathbb{C}}

\newtheorem{lem}{Lemme}
\newtheorem{theo}{Th\'eor\`eme}
\newtheorem{prop}{Proposition}
\newtheorem{cor}{Corollaire}[theo]

\begin{document}

\author{Florian Deloup}

\address{Universit\'e Paul Sabatier,
Toulouse III, Laboratoire \'Emile Picard de Math\'ematiques\\118,
route de Narbonne, 31062 Toulouse, France.}
\asciiaddress{Universite Paul Sabatier,
Toulouse III, Laboratoire Emile Picard de Mathematiques\\118,
route de Narbonne, 31062 Toulouse, France.}

\title{ Mono\"{\i}de des enlacements et facteurs orthogonaux}
\covertitle{Monoide des enlacements et facteurs orthogonaux}
\asciititle{Monoide des enlacements et facteurs orthogonaux 
(Monoids of linking pairings and orthogonal summands)}

\email{deloup@picard.ups-tlse.fr}

\begin{abstract}

A linking pairing is a symetric bilinear pairing $\lambda\co G \times
G \to \Q/\Z$ on a finite abelian group.  The set of isomorphism
classes of linking pairings is a non-cancellative monoid
$\mathfrak{E}$ under orthogonal sum, which is infinitely generated and
infinitely related. We propose a new presentation of $\mathfrak{E}$
that enables one to detect whether a linking pairing has a given
orthogonal summand. The same method extends to the monoid
${\mathfrak{Q}}$ of quadratic forms on finite abelian groups. We
obtain a combinatorial classification of ${\mathfrak{Q}}$ (that was
previously known for groups of period $4$).

As applications, we describe explicitly $3$-manifolds having a
degree one map onto prescribed (or proscribed) lens spaces. Most
of the results extend to $3$-manifolds endowed with a
parallelization or a spin structure. In particular, the
Reidemeister--Turaev function detects the existence of a spin
preserving degree one map between a rational homology $3$-sphere
and a lens space.

{\bf R\'esum\'e}

Un enlacement est une forme bilin\'eaire sym\'etrique $\lambda\co
G \times G \to \Q/\Z$ sur un groupe ab\'elien fini. L'ensemble des
classes d'isomorphismes d'enlacements forme un mono\"{\i}de
$\mathfrak{E}$, pour la somme orthogonale, \`a un nombre infini de
g\'en\'erateurs et de relations, sans simplification. Nous
proposons une nouvelle pr\'esentation de ${\mathfrak{E}}$ qui
permet de reconna\^{\i}tre si un enlacement poss\`ede un facteur
orthogonal donn\'e. La m\^eme m\'ethode se g\'en\'eralise au
mono\"{\i}de ${\mathfrak{Q}}$ des formes quadratiques sur les
groupes ab\'eliens finis. Nous obtenons ainsi une classification
combinatoire de ${\mathfrak{Q}}$, classification qui n'\'etait
pr\'ec\'edemment connue que pour les groupes de p\'eriode $4$.

Comme application, nous d\'ecrivons explicitement les
3-vari\'et\'es admettant une application de degr\'e un sur des
lenticulaires prescrits (ou proscrits). La plupart des r\'esultats
se g\'en\'eralisent aux 3-vari\'et\'es munies d'une
parall\'elisation ou d'une structure spinorielle. En particulier,
la fonction de Reidemeister--Turaev distingue l'existence ou non
d'une application de degr\'e un pr\'eservant les structures
spinorielles entre une 3-sph\`ere d'homologie rationnelle et un
lenticulaire.

\end{abstract}

\asciiabstract{%
A linking pairing is a symetric bilinear pairing lambda: GxG --> Q/Z
on a finite abelian group.  The set of isomorphism classes of linking
pairings is a non-cancellative monoid E under orthogonal sum, which is
infinitely generated and infinitely related. We propose a new
presentation of E that enables one to detect whether a linking pairing has
a given orthogonal summand. The same method extends to the monoid Q of
quadratic forms on finite abelian groups. We obtain a combinatorial
classification of Q (that was previously known for groups of period
4).  As applications, we describe explicitly 3-manifolds having a
degree one map onto prescribed (or proscribed) lens spaces. Most of
the results extend to 3-manifolds endowed with a parallelization or a
spin structure. In particular, the Reidemeister--Turaev function
detects the existence of a spin preserving degree one map between a
rational homology 3-sphere and a lens space.  }



 \maketitle

\section{Introduction}

Un enlacement $(G,b)$ est une forme bilin\'eaire sym\'etrique
$b\co G \times G \to \Q/\Z$ sur un groupe ab\'elien fini dont
l'homomorphisme adjoint $\widehat{b}\co G \to \hbox{Hom}(G,\Q/\Z)$
est bijectif. Les enlacements apparurent tout d'abord en topologie
comme invariants alg\'ebriques mesurant l'enlacement des
$(2n-1)$-cycles dans une vari\'et\'e ferm\'ee orient\'ee de
dimension $(4n-1)$. Une telle vari\'et\'e $M$ poss\`ede en effet
un enlacement $\lambda_{M}\co \hbox{Tors}\ H_{2n-1}(M) \times
\hbox{Tors}\ H_{2n-1}(M) \to \Q/\Z$ (voir \cite{Seifert}). H.
Minkowski avait d\'ej\`a indiqu\'e comment obtenir un syst\`eme
complet d'invariants de la classe d'isomorphismes de $\lambda_{M}$
\`a l'aide de syst\`emes d'\'equations de congruence et de sommes
de Gauss. Un syst\`eme complet d'invariants num\'e\-riques \'etait
connu de H. Seifert, dans le cas de $p$-groupes avec $p$ impair,
et de E. Burger \cite{Burger} dans le cas g\'en\'eral. La
classification compl\`ete fut ensuite poursuivie par C.T.C. Wall
\cite{Wall} dans le cas des $p$-groupes avec $p$ impair, puis
compl\'et\'ee par A. Kawauchi et S. Kojima \cite{KK} dans le cas
g\'en\'eral.

Deux enlacements $(G,b)$ et $(G',b')$ sont isomorphes s'il existe
un isomorphisme $\phi\co G \to G'$ tel que $b'(\phi(x), \phi(y)) =
b(x,y)$ pour tout $x, y \in G$. \'Etant donn\'es deux enlacements
$(G,b), (G',b')$, leur somme orthogonale $(G,b) \oplus (G',b')$,
not\'ee \'egalement $b \oplus b'$, est d\'efinie par $(b \oplus
b')((x,x'),(y,y')) = b(x,y) + b'(x',y')$ pour tous $x,y \in G$ et
$x', y' \in G'$. L'ensemble des classes d'isomorphismes
d'enlacements forme un mono\"{\i}de ${\mathfrak{E}}$ pour la somme
orthogonale $\oplus$. Ce mono\"{\i}de
 ${\mathfrak{E}}$ poss\`ede une infinit\'e de
g\'en\'erateurs et de relations et est sans simplification.

Une pr\'esentation par g\'en\'erateurs et relations de
${\mathfrak{E}}$ est propos\'ee dans \cite{KK}. La difficult\'e
majeure r\'eside dans les enlacements sur les $2$-groupes qui
forment un sous-mono\"{\i}de ${\mathfrak{M}} \subset
{\mathfrak{E}}$. Au-del\`a du caract\`ere r\'esiduellement fini de
${\mathfrak{E}}$, il ne semble pas exister de classification des
mono\"{\i}des dans laquelle ${\mathfrak{E}}$ s'ins\`ererait
naturellement. Nous proposons une autre approche pour \'eclaircir
la nature de~${\mathfrak{E}}$.

La motivation pour cette approche est topologique. En effet, le
syst\`eme complet d'invariants propos\'e dans \cite{KK} se
d\'eduit de th\'eories topologiques des champs ab\'eliennes
\cite{DG}\cite{De3}. Pour que ces th\'eories trouvent une
application proprement topologique, on cherche \`a d\'ecrire de
fa\c{c}on combinatoire comment reconstruire la classe
d'isomorphisme d'un enlacement \`a partir de ses invariants. Cette
id\'ee conduit \`a d\'ecrire l'``image" de ce syst\`eme
d'invariants. Comme cons\'equences, nous obtenons une nouvelle
pr\'esentation combinatoire de ${\mathfrak{M}}$, ainsi qu'un
algorithme pour reconna\^{\i}tre un facteur orthogonal d'un
enlacement. Nous g\'en\'e\-ralisons \'egalement cette approche au
cas des formes quadratiques (homog\`enes). Ceci permet d'obtenir
une pr\'esentation globale du mono\"{\i}de des formes quadratiques
sur les groupes finis, en particulier sur les $2$-groupes. Ceci
n'\'etait pr\'ec\'edemment connu que pour les $2$-groupes de
p\'eriode au plus $4$.

Nous appliquons cette m\'ethode pour d\'eterminer explicitement
des $3$-vari\'et\'es admettant des applications de degr\'e $1$ sur
des lenticulaires prescrits ou proscrits. Finalement, la plupart
des r\'esultats obtenus se g\'en\'eralise au cadre spinoriel,
c'est-\`a-dire aux $3$-vari\'et\'es munies d'une structure spin ou
d'une paral\-l\'elisation. En particulier, nous montrons que la
fonction de Reidemeister--Turaev \cite{Turaev} distingue
l'existence ou non d'une application de degr\'e un pr\'eservant
les structures spinorielles entre une 3-sph\`ere d'homologie
rationnelle et un lenticulaire. Certains des r\'esultats
alg\'ebriques pr\'esent\'es dans cet article ont
\'et\'e annonc\'es dans la note \cite{DeR}.

\textbf{Plan de l'article}\qua \S \ref{sec:comb} d\'ecrit une
pr\'esentation combinatoire du mono\"{\i}de ${\mathfrak{E}}$ des
enlacements, le cas le plus d\'elicat \'etant celui des
$2$-groupes (Th\'eor\`eme \ref{th:main}): ${\mathfrak{M}}$
s'identifie alors \`a un sous-mono\"{\i}de des fonctions $\N \to
(\Z/8\Z \cup \{ \infty \}) \times \N$, appel\'e le mono\"{\i}de
des tableaux admissibles. Cette description est appliqu\'ee dans
\S \ref{sec:reconnaissance} \`a la reconnaissance de facteurs
orthogonaux dans un enlacement (Th\'eor\`eme \ref{th:facteur}). \S
\ref{sec:quadratique} g\'en\'eralise les r\'esultats des sections
pr\'ec\'edentes  aux formes quadratiques (Th\'eor\`eme
\ref{th:classification-formesquad}). \S \ref{sec:applications}
pr\'esente quelques-unes des applications topologiques des
sections pr\'ec\'edentes dans les cas non-parall\'elis\'e et
parall\'elis\'e (ou spin). \S \ref{sec:preuve-main} contient les
d\'emonstrations des Th\'eor\`emes \ref{th:main} et
\ref{th:lenticulaire-proscrit}. Enfin, \S \ref{sec:questions}
pr\'esente quelques questions ouvertes.

Ce travail est en partie financ\'e par un contrat de recherche
Marie Curie de l'Union europ\'eenne MERG-CT-2004-510590.

\section{Une pr\'esentation combinatoire du mono\"{\i}de des
enlacements} \label{sec:comb}

\subsection{Le syst\`eme d'invariants de Minkowski-Burger}
\label{subsec:minkowski-burger}

Rappelons que tout enlacement $(G,\lambda)$ admet une unique d\'ecomposition
orthogonale
$$(G,\lambda) = \displaystyle \bigoplus_{p\ {\rm{prime}}} (G^{p}, \lambda^{p})$$
o\`u $(G^{p}, \lambda^{p})$ est un enlacement sur un $p$-groupe
(de type fini). Les cas $p >2$ et $p = 2$ sont distincts. Tout enlacement
$(G, \lambda)$ sur un $p$-groupe admet \`a son tour une
d\'ecomposition orthogonale
\begin{equation}
(G,\lambda) = \displaystyle \bigoplus_{k \geq 1} (G_{k}, \lambda_{k})
\label{eq:decompose}
\end{equation}
o\`u $(G_{k}, \lambda_{k})$ est un enlacement sur un
$\Z/p^{k}\Z$-module libre. Si $p \not= 2$, la d\'ecomposition est
unique (\`a isomorphisme pr\`es des facteurs
$(G_{k},\lambda_{k})$) et un tel enlacement est toujours isomorphe
\`a
 une somme orthogonale d'enlacements
cycliques sur des copies de $\Z/p^{k}\Z$. Le rang
$\rho_{k}(\lambda)$ de $G_{k}$ (en tant que $\Z/p^{k}\Z$-module)
est un invariant de $(G,\lambda)$; il est additif sur la somme
orthogonale. Au moyen de l'injection $1 \mapsto \frac{1}{p^{k}}$,
$\Z/p^{k}\Z \to \Q/\Z$, on peut regarder $\lambda_k$ comme un
enlacement \`a valeurs dans $\Z/p^{k}\Z$ (au lieu de $\Q/\Z$).
Ainsi le d\'eterminant $\det \lambda_k \in (\Z/p^{k}\Z)^{\times}$
est un \'el\'ement inversible bien d\'efini. On d\'efinit un
second invariant $\sigma_{k}(\lambda) \in \{ -1, 1 \} = \Z/2\Z$
comme le r\'esidu quadratique de $\det \lambda_k$ \penalty-2000
modulo $p^{k}$ (ou, ce qui revient au m\^eme, modulo $p$). Si
$\rho_{k}(\lambda) = 0$ alors $\lambda_{k} = 0$ et on convient que
$\sigma_{k}(\lambda) = 1$ ($0$ est un r\'esidu quadratique modulo
$p$). Regroupons alors les invariants ci-dessus sous la forme
d'une seule application $(\rho, \sigma)\co \N^{\times} \to \N
\times \Z/2\Z, k \mapsto (\rho_{k}(\lambda),
\sigma_{k}(\lambda))$. Le r\'esultat principal pour $p$ premier
impair est d\^u \`a Minkowski et sous la forme ci-dessous, \`a E.
Seifert et \`a C.T.C. Wall.

\begin{prop} \label{prop:injectivite}
Soit $p$ un nombre premier impair et
$(G,\lambda)$ un enlacement sur un $p$-groupe fini.
Le syst\`eme d'invariants ${\cal{S}} = (\rho(\lambda),\sigma(\lambda))$
d\'etermine la classe d'isomorphisme de $(G,\lambda)$.
\end{prop}

De plus, le syst\`eme est minimal en ce sens qu'\'etant donn\'ee toute
 sous-famille stricte ${\cal{F}} \subset {\cal{S}}$ d'invariants,
il existe des enlacements non distingu\'es par ${\cal{F}}$ qui sont
non isomorphes.

Le cas $p = 2$ est plus compliqu\'e, du fait que la
d\'ecomposition orthogonale (\ref{eq:decompose}) n'est pas unique
en g\'en\'eral. L'entier $\rho_{k}(\lambda)$ reste bien s\^ur un
invariant additif de l'enlacement. Un second invariant est
d\'efini \`a partir de sommes de Gauss. Soit $k \geq 1$.
Consid\'erons le nombre complexe
$$ \Gamma_{k}(G,\lambda) = \sum_{x \in G} \exp(\sqrt{-1}\ \pi\ 2^{k} \lambda(x,x)).$$
Il est bien connu que si $\Gamma_{k}(G,\lambda)\not=0$ alors
$\frac{\Gamma_{k}(G,\lambda)}{|\Gamma_{k}(G,\lambda)|}$ est une
racine $8$-\`eme de l'unit\'e \cite[\S 2]{Scharlau}. On d\'efinit
alors
$$\sigma_{k}(\lambda) = \left\{ \begin{array}{cl}
\frac{1}{2\pi}{\rm{Arg}}\ \Gamma_{k}(G,\lambda) \in \Z/8\Z & {\rm{si}}\
\Gamma_{k}(G,\lambda) \not= 0 \\
\infty & {\rm{si}}\ \Gamma_{k}(G,\lambda) = 0
\end{array} \right.$$
Soit $\overline{\Z}_{8} = \Z/8\Z \cup \{ \infty \}.$ Il s'agit du
mono\"{\i}de obtenu en adjoignant au groupe cyclique \`a 8
\'el\'ements un \'el\'ement suppl\'ementaire not\'e $\infty$, avec
la r\`egle $\infty + a = a + \infty = \infty = \infty + \infty$
pour tout $a \in \Z/8\Z$. Comme les sommes de Gauss ci-dessus sont
multiplicatives sur les sommes orthogonales, chaque $\sigma_k$
d\'efinit un homomorphisme de mono\"{\i}des ${\mathfrak{M}} \to
\overline{\Z}_{8}$. Chaque enlacement $\lambda$ donne ainsi lieu
\`a  une application $\sigma(\lambda)\co \N^{\times} \to
\overline{\Z}_{8}, \ k \mapsto \sigma_{k}(\lambda)$. On peut ainsi
\`a nouveau regrouper les invariants $\rho_k$ et $\sigma_k$ sous
la forme d'une seule application $(\rho, \sigma)\co \N^{\times}
\to \N \times \overline{\Z}_{8}$. Le r\'esultat principal de
classification par invariants est d\^u \`a E. Burger et sous la
forme ci-dessous, \`a A. Kawauchi et S. Kojima \cite[Th\'eor\`eme
4.1]{KK}.

\begin{prop} \label{prop:injectivite-2}
Soit $(G,\lambda)$ un enlacement sur un $2$-groupe fini.
Le syst\`eme d'invariants ${\cal{S}} = (\rho(\lambda), \sigma(\lambda))$
d\'etermine la classe d'isomorphisme de $(G,\lambda)$.
\end{prop}

L\`a encore, il est ais\'e de se rendre compte
 que le syst\`eme ${\cal{S}}$ est minimal.

\subsection{Le mono\"{\i}de des enlacements} \label{sec:monoide}

Soit ${\cal{M}}$ un mono\"{\i}de additif et $I$ un suite d'entiers
cons\'ecutifs.
Un tableau est une application $I \to {\cal{M}}$, qu'il sera pratique de
consid\'erer comme un diagramme de la forme
$$ \begin{array}{|c|c|c|c|c|} \hline
  I & k & k+1 & \ldots & l  \\ \hline \hline
{\cal{M}} & m_k & m_{k+1} & \ldots & m_{l}
\\
\hline
\end{array}$$
Afin de simplifier la notation, les notations de l'intervalle
ainsi que du mono\"{\i}de seront omises des tableaux suivants. La
longueur d'un tableau $T$ est l'entier $1 + \sup_{(m,n) \in
I\times I} |m - n| \in \overline{\N} = \N \cup \{ \infty \}$. Un
tableau $T'$ est un prolongement d'un tableau $T$ si $T'$ prolonge
$T$ en tant qu'application. Dans ce cas, $T$ est un tableau
extrait de $T$. \'Etant donn\'e un tableau $T\co I \to {\cal{M}}$
quelconque, on peut toujours le prolonger trivialement sur $\N$
entier en d\'efinissant $\tilde{T}(n) = 0$ pour $n \in \N - I$. En
pratique, on confondra un tableau $T$ et son prolongement trivial
$\tilde{T}$ \`a $\N$ ainsi d\'efini. Ainsi on dira qu'un tableau
$T$ est {\it{fini}} s'il est de longueur finie ou s'il est le
prolongement trivial $\tilde{T}'$ d'un tableau $T'$ de longueur
finie. (C'est la d\'efinition habituelle de support fini.) Comme
${\cal{M}}$ est un mono\"{\i}de, l'addition de tableaux est bien
d\'efinie. La somme de deux tableaux $T_{1}$ et $T_{2}$ est
d\'efinie sur $\N$ par $T_{1} + T_{2} = \tilde{T}_{1} +
\tilde{T}_{2}$ o\`u $\tilde{T}_{i}$, $i = 1,2$, d\'esigne le
prolongement trivial \`a $\N$. L'ensemble des tableaux $\N \to
{\cal{M}}$ forme un mono\"{\i}de. L'\'el\'ement neutre $0$ est le
tableau envoyant $\N$ sur $0$. Le d\'elimiteur \`a gauche (resp.
\`a droite) d'un tableau $T\co  I \to {\cal{M}}$ est l'\'el\'ement
$-1 \leq {\rm{Inf}}\ I - 1 < \infty$ (resp. l'\'el\'ement $0 \leq
{\rm{Sup}}\ I + 1 \leq \infty$).

Soit $f\co \N^{\times} \to {\cal{M}}$ une application invariante
sur les classes d'isomorphismes d'enlacements. Nous dirons qu'un
tableau $T\co I \to {\cal{M}}$ est {\emph{admissible}} s'il existe
un enlacement $(G,\lambda)$ tel que $f(m) = T(m)$ pour tout $m \in
I$.

\subsubsection{Cas $p \not= 2$}

Nous commen\c{c}ons par le cas, techniquement plus simple, des
enlacements sur un $p$-groupe avec $p$ premier impair. Nous
consid\'erons le syst\`eme d'invariants $(\rho, \sigma)$; les
tableaux correspondants sont donc \`a valeurs dans le mono\"{\i}de
${\cal{M}} = \N \times \Z/2\Z = \N \times \{ \pm 1 \}$. Nous avons
vu (\S \ref{subsec:minkowski-burger}) que $\rho_{k}(\lambda) = 0$
implique $\sigma_k(\lambda) = 1$. Une condition n\'ecessaire
d'admissibilit\'e d'un tableau $T = (r(m), s(m))_{m \in
\N^{\times}}$ est donc
\begin{equation}
{\hbox{Pour tout}}\ m \in \N^{\times}, \ \ \ r(m) = 0
\Longrightarrow s(m) = 1. \label{eq:rang-signature}
\end{equation}

\begin{theo} \label{th:cas-impair-facile}
Soit $p$ premier distinct de $2$. Alors un tableau fini $T=(r,s)$
est admissible pour un enlacement sur un $p$-groupe si et
seulement si la condition $(\ref{eq:rang-signature})$ est
v\'erifi\'ee.
\end{theo}

Le Th\'eor\`eme \ref{th:cas-impair-facile} est bas\'e sur l'
observation suivante: si le rang $\rho_{k}(\lambda)$ est fix\'e,
alors la classe d'isomorphisme de $\lambda_{k}$ dans (\ref{eq:decompose})
d\'etermine et est d\'etermin\'ee par $\sigma_{k}(\lambda_{k}) \in \Z/2\Z$.

\subsubsection{Cas $p=2$}

Consid\'erons \`a pr\'esent les enlacement sur les $2$-groupes.
Les tableaux sont \`a valeurs dans le
 mono\"{\i}de ${\cal{M}} = \N \times \overline{\Z}_{8}$, que nous noterons
$T\co m \mapsto (r(m), s(m))$, avec $r(m) \in \N$ (rang formel)
and $s(m) \in \overline{\Z}_{8}$ (signature formelle).
 Nous dirons qu'un tableau est {\it{admissible}} s'il existe un enlacement
$(G,\lambda)$ sur un $2$-groupe tel que $r(m) = \rho_{m}(\lambda)$ et $s(m)
 = \sigma_{m}(\lambda)$ pour tout
$m \in I$.

Un entier $m \in I$ sera dit {\emph{r\'egulier}} pour un tableau $T$ si
$r(m) = 0$ ou $s(m)\not= \infty$. On note $I_{\hbox{\scriptsize reg}}
\subseteq I$ l'ensemble des
\'el\'ements r\'eguliers de $T$.
Pr\'esentons quatre types particuliers distincts de tableaux:
\begin{enumerate}
\item[$\bullet$] Type T$_{0}$. Tout tableau de
longueur impaire de la forme $T=(0,s(m))_{m \in I}$.
\item[$\bullet$] Type T$_{1}$. Tout tableau de la forme
$\begin{array}{|c|} \hline m \\ \hline \hline
1 \\ \hline
\infty \\ \hline
 \end{array}$ pour un entier non nul $m$.
\item[$\bullet$] Type T$_{2}$. Tout tableau de la forme
$\begin{array}{|c|} \hline m \\ \hline \hline
2 \\ \hline
\infty \\ \hline
 \end{array}$ pour un entier non nul $m$.
\item[$\bullet$] Type T$_{3}$. Tout tableau
de longueur impaire tel que $I = I_{\hbox{\scriptsize reg}}$.
\end{enumerate}

Le r\'esultat principal est un crit\`ere n\'ecessaire et suffisant pour qu'un tableau soit
admissible.

\begin{theo} \label{th:main}
Un tableau fini $T\co \N^{\times} \to \N \times
\overline{\Z}_{8},\ m \mapsto (r(m), s(m))$ est admissible si et
seulement si les conditions suivantes sont satisfaites:
\begin{enumerate}
\item[$(1)$] $r(I_{\rm{{reg}}}) \subseteq 2\N$.
\item[$(2)$] $s(m) = \sum_{k \geq m+1} r(k)\ ({\rm{mod}}\ 2)$ pour tout $m \in I_{\rm{ reg}}$.
\item[$(3)$] $s(m) + s(m+1) =
2 \sum_{k \geq m+2} r(k)\ ({\rm{mod}}\ 4)$ pour tout $\{m, m+1 \} \subseteq I_{\rm{ reg}}$.
\item[$(4)$] Pour tout tableau $T_{\rm{ ext}}$ extrait de $T$ et
 pour toute paire de
d\'elimiteurs $m, n$ de $T_{\rm{ ext}}$ dans
$I_{\rm{{reg}}}$,
les conditions
suivantes sont v\'erifi\'ees:
$$\begin{array}{|c||c|c|c|c|} \hline
{\rm{Type\ de}}\ T_{\rm{ ext}} & {\rm{T}}_{0} & {\rm{T}}_{1} & {\rm{T}}_{2} & {\rm{T}}_{3} \\ \hline
s(m) - s(n) & 0 & \pm 1 & 0, \pm 2 & 0, 4\\ \hline
\end{array}$$
\end{enumerate}
\end{theo}

Compte-tenu du fait que le groupe d'un enlacement est fini, il est
ais\'e d'observer sur le rang et la signature que tout tableau
admissible est fini. Ceci garantit en particulier que les sommes
intervenant dans les conditions $(2)$ et $(3)$ sont finies. (En
particulier, la condition $(2)$ implique que $s(m)\not= \infty$
d\`es que
 $r(m) =
0$: les entiers r\'eguliers $m$ de $T$ sont exactement les entiers $m$ tels
que $s(m) \not= \infty$.)
De mani\`ere g\'en\'erale, la n\'ecessit\'e des conditions
\'enonc\'ees dans le Th\'eor\`eme \ref{th:main} est une cons\'equence de
calculs classiques d'enlacements
et de sommes de Gauss. La preuve de la suffisance est constructive et
sera donn\'ee en \S \ref{sec:preuve-main}.

Notons ${\mathfrak{T}}$ le mono\"{\i}de constitu\'e des tableaux
$T\co \N^{\times} \to \N \times \overline{\Z}_{8}$. On d\'eduit du
Th\'eor\`eme \ref{th:main} que la somme de deux tableaux
admissibles est encore un tableau admissible, de sorte que le
sous-ensemble des tableaux admissibles constitue un
sous-mono\"{\i}de ${\mathfrak{T}}^{\rm{ adm}}$ de
${\mathfrak{T}}$. Puisque $\rho,\sigma$ sont des invariants
complets du mono\"{\i}de ${\mathfrak{M}}$ des classes
d'isomorphismes d'enlacements sur les 2-groupes, l'application
$(\rho, \sigma)\co {\mathfrak{M}} \to {\mathfrak{T}}$ est
injective. Il en r\'esulte la description combinatoire de
${\mathfrak{M}}$ ci-dessous.

\begin{cor} \label{th:prescomb}
Le mono\"{\i}de ${\mathfrak{M}}$  des classes d'isomorphismes d'enlacements sur les 2-groupes
est isomorphe au sous-mono\"{\i}de ${\mathfrak{T}}^{\rm{ adm}}$ des tableaux
admissibles.
\end{cor}

Les Th\'eor\`emes \ref{th:cas-impair-facile} et \ref{th:main}
ensemble donnent ainsi une pr\'esentation combinatoire compl\`ete
du mono\"{\i}de ${\mathfrak{E}}$ des enlacements.

Le Th\'eor\`eme \ref{th:main} permet de calculer le nombre de
classes d'isomorphismes d'enlace\-ment ayant un rang ou une
signature donn\'e, tout au moins th\'eoriquement. Je ne connais
pas de formule explicite. Les quelques remarques suivantes peuvent
\^etre utiles. Tout d'abord, on peut d\'efinir deux applications
``profil'' par
profil$_{\rho}(\lambda)=(k,\rho_{k}(\lambda))_{k\geq 1}$ et
profil$_{\sigma}(\lambda)=(k,\sigma_{k}(\lambda))_{k \geq 1}$ et
\'etudier les fibres de ces applications. Y a-t-il en particulier
des fibres ``g\'en\'eriques'' ?  L'approche la plus encourageante
semble \^etre l'\'etude des fibres de profil$_{\rho}$. Plus
globalement, d\'efinissons alors l'application Profil par
$${\rm{Profil}}(\lambda) = \{ k \in \N\ | \ \rho_{k}(\lambda)
\not= 0 \}.$$ Est-il possible de classifier ${\mathfrak{M}}$ \`a
partir des fibres de Profil  ? Dans ce contexte, les lemmes 3.1,
3.2 et 3.3 de \cite{KK} s'interpr\`etent comme la classification
des fibres Profil$^{-1}(\{k \})$, $k \geq 1$ et la proposition 5.2
de \cite{KK}
comme le calcul du groupe de Witt $W({\rm{Profil}}^{-1}(\{k \}))$.

On a d\'ej\`a observ\'e que la d\'ecomposition orthogonale d'un enlacement sur un $2$-groupe
n'est pas unique, m\^eme \`a isomorphisme pr\`es des facteurs (et m\^eme pour un $2$-groupe
homog\`ene, c'est-\`a-dire isomorphe \`a une somme directe de copies d'un groupe cyclique).
On peut cependant montrer que tout enlacement sur un $2$-groupe homog\`ene admet une forme
``normale'' qui est unique: voir \cite[\S 3]{KK} et \cite[\S 3]{Mi}. R. Miranda a en fait montr\'e
qu'il existe une forme normale pour un enlacement sur un $2$-groupe quelconque \cite[\S 4]{Mi}.

\section{La reconnaissance d'un facteur orthogonal dans un enlacement}
\label{sec:reconnaissance}

Consid\'erons \`a pr\'esent la question de reconna\^{\i}tre si un
enlacement $\lambda'$ est un
facteur orthogonal d'un enlacement $\lambda$, c'est-\`a-dire s'il
existe un enlacement $\lambda"$
tel que $$\lambda = \lambda' \oplus \lambda".$$ D\'ecrivons tout d'abord
des conditions n\'ecessaires
simples pour qu'une telle d\'ecomposition orthogonale existe. Il est
clairement n\'ecessaire que
\begin{equation}
\rho_{k}(\lambda) \geq \rho_{k}(\lambda')\ {\hbox{pour tout}}\ k \geq 1. \label{eq:rankgeq}
\end{equation}
Une seconde condition n\'ecessaire r\'esulte du comportement
 de $\sigma$ sur les
sommes orthogonales. Dans le cas $p \not= 2$, au vu du
Th\'eor\`eme \ref{th:cas-impair-facile}, il existe toujours un
enlacement sur un $p$-groupe de signature formelle prescrite
pourvu que le rang formel soit non nul. On en d\'eduit:

\begin{theo}
Soit $(G,\lambda)$ un enlacement sur un $p$-groupe fini.
L'enlacement $\lambda'$ est un facteur orthogonal de $\lambda$ si et
seulement si
$${\rm{pour}}\ {\rm{tout}}\ k \geq 1, \ \ \ \left\{  \begin{array}{c}
\rho_{k}(\lambda') < \rho_{k}(\lambda),\ {\rm{ou}}\\
\rho_{k}(\lambda') = \rho_{k}(\lambda)\ {\rm{and}}\ \sigma_{k}(\lambda') =
\sigma_{k}(\lambda).
\end{array} \right. $$
\end{theo}

Dans le cas $p = 2$, l'additivit\'e de $\sigma$ implique
\begin{equation}
\sigma_{k}(\lambda') = \infty \ \ \Longrightarrow \ \ \sigma_{k}(\lambda) = \infty,\ \ \  {\hbox{pour tout}}\ k \geq 1. \label{eq:inftyinfty}
\end{equation}
Supposons \`a pr\'esent ces conditions (\ref{eq:rankgeq}) et
(\ref{eq:inftyinfty}) v\'erifi\'ees. Soit $E(\lambda')$ l'ensemble
des applications de $\{ m \in \N^{\times} \ | \
\sigma_{m}(\lambda') = \infty \}$ dans $\overline{\Z}_{8}$. Nous
allons associer \`a $(\lambda, \lambda')$ un ensemble
$$S_{\lambda, \lambda'} = \{ T_{\alpha} \}_{\alpha \in
\overline{\Z}_{8}}$$ de tableaux. Pour $a \in E(\lambda')$, nous
d\'efinissons le tableau $T_{a} = (r_{a}, s_{a})\co \N^{\times}
\to \N \times \overline{\Z}_{8}$ par
\begin{equation}
\begin{array}{lcl}
r_{a}(k) & = & \rho_{k}(\lambda) - \rho_{k}(\lambda')\\
s_{a}(k)& = & \left\{
 \begin{array}{cl}
a(k) & \hbox{si}\ \sigma_{k}(\lambda') = \infty \\
 \sigma_{k}(\lambda) - \sigma_{k}(\lambda') & \hbox{si}\ \sigma_{k}(\lambda') \not=
\infty
  \end{array}
\right.
\end{array} \ \ \ {\hbox{pour tout $k \in \N^{\times}$.}}
\end{equation}
Le tableau $T_{a}$ est bien d\'efini gr\^ace \`a la condition
(\ref{eq:rankgeq}) et au fait que $\infty$ est le seul \'el\'ement
non inversible dans $\overline{\Z}_8$.

\begin{theo} \label{th:facteur}
Un enlacement $\lambda'$ est un facteur orthogonal d'un enlacement $\lambda$ si et
seulement si les conditions $(\ref{eq:rankgeq})$ et $(\ref{eq:inftyinfty})$ ci-dessus
sont v\'erifi\'ees et
s'il existe un tableau admissible $T \in S_{\lambda, \lambda'}$.
\end{theo}

\proof[D\'emonstration] Si $\lambda = \lambda' \oplus
\lambda''$, on v\'erifie que le tableau $$T_{\lambda''} =
(\rho_{k}(\lambda''), \sigma_{k}(\lambda''))_{k \in
\N^{\times}}$$ d'invariants associ\'e \`a $\lambda''$ est dans
$S_{\lambda, \lambda'}$. R\'eciproquement, si $T$ est admissible,
d'apr\`es le Th\'eor\`eme \ref{th:main}, il existe un enlacement
$\lambda''$ dont $T = T_{\lambda''} = (\rho_{k}(\lambda''),
\sigma_{k}(\lambda''))_{k \in \N^{\times}}$ est le tableau des
invariants. On v\'erifie imm\'ediatement la relation suivante, au
niveau des tableaux d'invariants, respectivement de $\lambda,
\lambda'$ et $\lambda''$:
$$T_{\lambda} = T_{\lambda'} + T_{\lambda''} = T_{\lambda' \oplus \lambda''},$$
o\`u la derni\`ere \'egalit\'e r\'esulte de l'additivit\'e des
invariants $\rho$ et $\sigma$ sur les sommes orthogonales.
L'application qui \`a (une classe d'isomorphisme d') un enlacement
associe ses invariants $(\rho, \sigma)$ \'etant injective
(Proposition \ref{prop:injectivite-2}), on conclut que
$\lambda = \lambda' \oplus \lambda''$. \endproof

\section{Le mono\"{\i}de des formes quadratiques} \label{sec:quadratique}

 Consid\'erons bri\`evement le cas
plus g\'en\'eral
 des formes quadratiques
sur un groupe ab\'elien fini. Une forme quadratique sur un groupe
ab\'elien fini $G$ est une application $q\co G \to \Q/\Z$ telle
que $q(nx) = n^{2}q(x)$ pour tout $(n, x) \in \Z \times G$ et
telle que l'application $\lambda_{q}\co G \times G \to \Q/\Z$
d\'efinie par $\lambda_{q}(x,y) = q(x+y) - q(x) - q(y)$ soit un
enlacement. Les formes quadratiques $G \to \Q/\Z$ ayant le m\^eme
enlacement associ\'e sont en bijection avec
${\rm{Hom}}(G,\Z/2\Z)$. Il en r\'esulte que sur le facteur
orthogonal $G_{\rm{impair}}$ des \'el\'ements d'ordre impair, les
formes quadratiques sont d\'etermin\'ees par leur enlacement
associ\'e. Consid\'erons alors les formes quadratiques sur les
$2$-groupes. Il r\'esulte de \cite[Th.~5]{Wall} qu'une telle forme
quadratique $q\co G \to \Q/\Z$ est classifi\'ee par les invariants
$\rho_{k}(\lambda_{q}), \sigma_{k}(\lambda_{q})$ associ\'ees \`a
l'enlacement associ\'e $\lambda_{q}$ et un seul invariant
suppl\'ementaire, la somme de Gauss $\gamma(q) = \sum_{x \in G}
\exp(2 i \pi q(x)) \in \C$. Aussi la construction combinatoire \`a
l'aide des tableaux est essentiellement la m\^eme: on consid\`ere
maintenant le mono\"{\i}de ${\mathfrak{T}}$ constitu\'e des
tableaux $(r,s)\co \N \to \N \times \overline{\Z}_{8}$. Le tableau
$T_{q} = (\rho,\sigma)\co \N \to
 \N \times \overline{\Z}_{8}$ d'invariants associ\'e \`a $q$ est d\'efini par
$$\rho_{k}(q) = \rho_{k}(\lambda_{q}) \ \hbox{pour}\ k \geq 1
\ \ \hbox{et}\ \ \rho_{0}(q) = 0$$ $$\sigma_0(q) =
\frac{1}{2\pi}{\rm{Arg}}\left( \gamma(q) \right) \in \Z/8\Z \
\hbox{et}\ \ \sigma_{k}(q) = \sigma_{k}(\lambda_{q})\in
\overline{\Z}_{8}\ \hbox{pour}\ k \geq 1.\leqno{\rm et}$$ (Noter que comme
$\lambda_{q}$ est non d\'eg\'en\'er\'ee, $\gamma(q) \not= 0.$)
Avec cette modification, le Th\'eor\`eme \ref{th:main}, le
corollaire \ref{th:prescomb} ainsi que le Th\'eor\`eme \ref{th:facteur}
s'\'etendent au cas quadratique. Nous obtenons
en particulier le th\'eor\`eme suivant.

\begin{theo} \label{th:classification-formesquad}
Le mono\"{\i}de des formes quadratiques sur les $2$-groupes finis est
isomorphe au sous-mono\"{\i}de constitu\'e des tableaux $$\N \to \N \times
\overline{\Z}_{8},\ m \mapsto (r(m), s(m))$$ v\'erifiant $r(0) = 0$ et
$s(0) \in \Z/8\Z$ ainsi que les conditions $(1)$ \`a $(4)$ du Th\'eor\`eme
\ref{th:main}.
\end{theo}

\textbf{Remarque}\qua Ce r\'esultat donne une
pr\'esentation globale du mono\"{\i}de des formes quadratiques sur les
groupes finis.
En particulier, le Th\'eor\`eme \ref{th:classification-formesquad}
g\'en\'eralise les pr\'esentations connues du mono\"{\i}de des
formes quadratiques sur les $2$-groupes de p\'eriode $2$ ou $4$,
voir \cite[sec.\ 3.4.3, Th.\ 3.6.5]{DIK}. Le mono\"{\i}de des formes
quadratiques sur les $2$-groupes de p\'eriode $2$ permet de classifier
les surfaces immerg\'ees \`a homotopie r\'eguli\`ere pr\`es \cite[Th.~4]{PIN}.

Si l'on appelle un tableau admissible un tableau v\'erifiant les
conditions du Th\'eor\`eme \ref{th:classification-formesquad}, le
Th\'eor\`eme \ref{th:facteur} plus haut se g\'en\'eralise mutatis
mutandis aux formes quadratiques.

\section{Quelques applications} \label{sec:applications}

\'Etant donn\'ee une vari\'et\'e orient\'ee $M$ de dimension
$4n-1$, on note $\lambda_{M}$ son enlacement sur ${\rm{Tors}}\ H^{2n}(M)$.

\subsection{Lenticulaires et facteurs d'enlacements}

\begin{prop} \label{prop:ortho}
Soit $f\co M \to X$ une application de degr\'e $d$ entre deux
vari\'et\'es diff\'erentiables ferm\'ees orient\'ees connexes de
dimension $4n-1$. Alors $d \cdot {\rm{Ker}}\ \left(f^{*}\co
{\rm{Tors}}\ H^{2n}(X) \to {\rm{Tors}}\ H^{2n}(M) \right) = 0$. En
particulier, si $d$ est premier avec l'exposant de ${\rm{Tors}}\
H^{2n}(M)$, alors l'enlacement $\lambda_{X}$ est un facteur
orthogonal de $\lambda_{M}$.
\end{prop}

\proof[D\'emonstration] La naturalit\'e en cohomologie fournit
la relation $\lambda_{M} \circ (f^{*} \times f^{*}) = d \ \lambda_{X}$.
Ceci implique la premi\`ere affirmation. Si $d$ est premier
avec l'exposant de ${\rm{Tors}}\ H^{2n}(M)$, alors l'application $f^*$
est injective. Donc $f^*({\rm{Tors}}$ $H^{2n}(X))$ est un sous-groupe
de ${\rm{Tors}}\ H^{2n}(M)$ sur la restriction duquel $\lambda_{M}$ est
non-singulier. Le r\'esultat s'ensuit. \endproof

Sur les $3$-vari\'et\'es elliptiques, on peut montrer une r\'eciproque.
En particulier, on a le r\'esultat suivant \cite{HWZ}.

\begin{theo} \label{th:hwz}
Il existe une application $f\co M \to L(n,p)$ de degr\'e $1$ d'une
$3$-vari\'et\'e ferm\'ee orient\'ee sur un espace lenticulaire si
et seulement si $\lambda_{M}$ contient l'enlacement de $L(n,p)$
comme facteur orthogonal. En particulier, l'un des facteurs
orthogonaux de $\lambda_{M}$ est cyclique.
\end{theo}

Nous allons utiliser ce dernier r\'esultat pour d\'ecrire certaines
$3$-vari\'et\'es admettant (resp. n'admettant pas) des applications de
degr\'e $1$ sur des
lenticulaires prescrits (resp. proscrits).

Le premier r\'esultat est une g\'en\'eralisation de
\cite[Prop. 6.1]{KK}, simple cons\'equence du Th\'eor\`eme  \ref{th:hwz}.

\begin{prop}
S'il existe une application $L(n,m) \# L(n,m') \to X^{3}$ de degr\'e $1$ o\`u
$X$ se plonge de fa\c{c}on lisse dans $S^{4}$, alors $n$ est impair et
$L(n,m)$ et $-L(n,m')$ ont le m\^eme type d'homotopie orient\'ee.
\end{prop}

On se propose maintenant de d\'eterminer \`a quelles conditions
une $3$-vari\'et\'e $M$ ferm\'ee orient\'ee admet une application
de degr\'e $1$ sur {\emph{tout}} lenticulaire dont le groupe
fondamental $\pi$ est fix\'e. Il r\'esulte des consid\'erations
pr\'ec\'edentes qu'il suffit de consid\'erer le cas o\`u l'ordre
de  $\pi$ est une puissance d'un nombre premier $p$. Si $p$ est
impair, il faut et il suffit que $H_{1}(M)$ ait un facteur
orthogonal isomorphe \`a une somme d'au moins deux copies de
$\pi$. Dans le cas o\`u $p = 2$ et $\pi = \Z/p^{k}\Z$ avec $k \geq
3$, on peut donner une r\'eponse compl\`ete en utilisant le
Th\'eor\`eme \ref{th:facteur}. \'Etant donn\'e un intervalle fini
$I$ de $\N$ et $a \in I$, le sym\'etrique $I'$ de $I$ est d\'efini
par $a+k \in I'$ si et seulement si $a - k \in I$. \'Etant donn\'e
un tableau fini $T:I \to \N \times \overline{\Z}_{8}$ et $a \in
I$, son {\emph{sym\'etrique}} par rapport \`a $a$ est le tableau
$T': I' \to \N \times \overline{\Z}_{8}$ o\`u $I'$ est le
sym\'etrique de $T$ par rapport \`a $a$ et $T'(a+k) = T(a-k)$ pour
tout $k \in \{ k \in \N \ |\ a+k \in  I'\}$.

Consid\'erons la liste ${\mathfrak{L}}$ ci-dessous de tableaux
$T\co I \to \N \times \overline{\Z}_{8},\ m \mapsto (r(m), s(m))$.
Le symbole $\N$ d\'esigne un entier positif ou nul arbitraire et
$\overline{\Z}_{8}$ un \'el\'ement arbitraire de
$\overline{\Z}_{8}$.
$$ \begin{array}{|c|} \hline
k \\
\hline \hline
r(k) \geq 4 \\
\hline
\infty \\
\hline
\end{array}\ ,\ \begin{array}{|c|c|c|}\hline
k-1 & k & k+1 \\
\hline \hline
0 & 3 & 0 \\
 \hline
s(k-1) & \infty & s(k+1) \\
\hline
\end{array}\ \ \ \hbox{\scriptsize ${\rm{avec}}\ s(k-1)-s(k+1) = \pm 1$,}$$
$$\begin{array}{|c|c|} \hline
k-1 & k \\
\hline \hline
r(k-1) \geq 1 & 3 \\
\hline
\overline{\Z}_8 & \infty \\
\hline
\end{array}\ ,\ \begin{array}{|c|c|c|} \hline
k-2 & k-1 & k \\
\hline \hline
r(k-2) \geq 1 & 0 & 3 \\
\hline
\overline{\Z}_8 & \Z_8 & \infty \\
\hline
\end{array}\ , \ \begin{array}{|c|c|} \hline
k-1 & k \\
\hline \hline
r(k-1) \geq 1 & 2 \\
\hline
 \infty & \infty \\
\hline
\end{array}\ ,
$$
$$ \begin{array}{|c|c|c|} \hline
k-1 & k & k+1 \\
\hline \hline
r(k-1) \geq 2 & 1 & r(k+1) \geq 1 \\
\hline
\overline{\Z}_8 & \infty & \infty \\
\hline
\end{array}\ , \ \begin{array}{|c|c|c|c|} \hline
k-2 & k-1 & k & k+1 \\
\hline \hline
r(k-2) \geq 1 & \N & 1 & r(k+1) \geq 1 \\
\hline
\infty &  \overline{\Z}_{8} & \infty & \infty \\
\hline
\end{array}.$$

\begin{theo}
Soit $k \geq 3$. Une $3$-vari\'et\'e $M$ ferm\'ee orient\'ee admet
une application de degr\'e $1$ sur {\emph{tout}} lenticulaire dont
le groupe fondamental $\pi$ est $\Z/2^{k}\Z$ si et seulement si
l'un des tableaux de la liste ${\mathfrak{L}}$, ou son
sym\'etrique par rapport \`a $k$, est un tableau extrait du
tableau $T = (\rho_{k}(\lambda_{M}), \sigma_{k}(\lambda_{M}))_{k
\in \N^{*}}$.
\end{theo}

\proof[D\'emonstration] V\'erifier que si $T$ prolonge l'un des
tableaux de $L$ alors $M$ admet une application de degr\'e $1$ sur tout
lenticulaire $L(2^{k},a)$ ne pose pas de probl\`eme particulier. Pour la
r\'eciproque, on utilise le Th\'eor\`eme \ref{th:facteur} en distinguant
les cas $\rho_k(\lambda_{M}) \geq 4$, $\rho_k(\lambda_M) =3$, $2$ ou $1$.
\endproof 

\textbf{Exemple}\qua Soit $a, b, c > 0$. Pour un entier $n
\geq 1$ et une vari\'et\'e $M$, on note $n \ M$ la somme connexe
$M \# \ldots \# M$ ($n$ fois). La $3$-vari\'et\'e $$a\
L(16,\alpha) \# b\ L(32,\beta) \# c\ L(64,\gamma)$$ admet une
application de degr\'e $1$ sur chaque lenticulaire $L$ tel que
$\pi_{1}(L) = \Z/32\Z$ si et seulement si l'une des conditions
suivantes est v\'erifi\'ee:
\begin{enumerate}
\item[$\bullet$] $b \geq 4$;
\item[$\bullet$] $b = 2$ et $a + c \geq 1$;
\item[$\bullet$] $b = 1$ et $a + c \geq 3$ et $ac \geq 2$.
\end{enumerate}

Dans une autre direction, nous avons le r\'esultat suivant.

\begin{theo}[Le lenticulaire proscrit]\label{th:lenticulaire-proscrit}
Soit $s$ un entier impair. Il existe une infinit\'e de
$3$-vari\'et\'es irr\'eductibles (hyperboliques) distinctes
admettant une application de degr\'e $1$ sur chaque lenticulaire
$L(16,r)$ pour $r \not\equiv s$ {\rm{mod}} $8$ et aucune
application de degr\'e $1$ sur $L(16,s)$.
\end{theo}

La d\'emonstration de ce dernier r\'esultat fait l'objet de la section \S
\ref{subsec:demth6}.

\subsection{Raffinements spinoriels et facteurs quadratiques}

% Notation multiplicative des actions....

Soit $M$ une $3$-vari\'et\'e ferm\'ee orient\'ee connexe. Il est connu que
le fibr\'e
tangent de $M$ est trivial.
 Une parall\'elisation de
$M$ est le choix d'une trivialisation $\tau$ de son fibr\'e tangent $T_{M}$
(consid\'er\'e \`a homotopie pr\`es). Le groupe $H^{1}(M;\Z/2\Z)$
agit librement et transitivement sur l'ensemble des
parall\'elisations de $M$. Dans la suite de ce paragraphe, les
groupes et les actions seront not\'es multiplicativement.
 Une structure spin sur $M$ est la donn\'ee
d'une trivialisation de $T_{M}$ sur son $1$-squelette qui
s'\'etend au $2$-squelette, consid\'er\'ee \`a homotopie pr\`es.
Il est clair que par restriction au $1$-squelette, une
trivialisation $t$ d\'etermine une structure spin.
R\'eciproquement, si une trivialisation s'\'etend au $2$-squelette
de $M$ alors elle s'\'etend en une trivialisation de $T_M$. \`A
toute structure spin $s$ de $M$, on sait associer de fa\c{c}on
canonique et naturelle une forme quadratique $q_{s}\co
{\rm{Tors}}\ H^{2}(M) \to \Q/\Z$ dont la forme bilin\'eaire
associ\'ee est $\lambda_{M}$ (voir \cite{LL} \cite{MS}). Le groupe
$H^{1}(M;\Z/2\Z)$ agit aussi sur les formes quadratiques
${\rm{Tors}}\ H^{2}(M) \to \Q/\Z$ via le Bockstein $\beta\co
H^{1}(M;\Z/2\Z) \to {\rm{Tors}}\ H^{2}(M)$. L'action explicite est
donn\'ee par la formule $(h \cdot q)(x) = q(x) + \lambda_M(\beta
h,x)$ pour tout $h \in H^{1}(M;\Z/2\Z)$ et $x \in {\rm{Tors}}\
H^{2}(M)$. On v\'erifie qu' elle est transitive. Elle est de plus
libre si $M$ est une sph\`ere d'homologie rationnelle.
 L'application $s \mapsto q_{s}$ ci-dessus est
$H^{1}(M;\Z/2\Z)$-\'equivariante
en ce sens que $$q_{h \cdot s} = h \cdot q_{s}, \ \ h \in H^{1}(M;\Z/2\Z).$$
Si $M$ est une $3$-sphere d'homologie rationnelle, alors
$s \mapsto q_s$ est bijective.

\begin{theo} \label{th:facteur-spin}
Soient $M, X$ deux $3$-sph\`eres d'homologie rationnelle
orient\'ees munies de structures spin $s_M$ et $s_X$
respectivement. On suppose que $X$ a le type d'homotopie d'un
lenticulaire. Alors il existe une application $f\co  M \to X$ de
degr\'e un telle que $f^{*}(s_{X}) = s_{M}$ si et seulement si
$q_{s_X}$ est un facteur orthogonal de $q_{s_M}$.
\end{theo}

\proof[D\'emonstration] Supposons l'existence de
l'application $f$ de degr\'e un comme dans l'\'enonc\'e. La
naturalit\'e de l'application $s \mapsto q_{s}$ fournit la
relation $q_{f^{*}s_X} = f^{*}q_{s_{X}} = q_{s_{X}} \circ f^{*} =
q_{s_{M}}$. On conclut alors par Prop. \ref{prop:ortho} (avec $d$
= 1). Pour la r\'eciproque,  par hypoth\`ese, on a une
d\'ecomposition orthogonale de la forme $q_{s_{M}} = q_{s_{X}}
\oplus q$. On peut construire un homomorphisme $\phi\co H^{2}(M)
\to H^{2}(X)$ tel que $q_{s_{X}} \circ \phi = q_{s_{M}}$ et induit
par un homomorphisme $\pi_{1}(M) \to \pi_{1}(X)$ (par dualit\'e de
Poincar\'e). Puisque $\pi_{2}(X) = 0$ et dim $X = 3$,
l'application naturelle $[M,X] \to {\rm{Hom}}(\pi_{1}(M),
\pi_{1}(X))$ est surjective (voir par exemple
\cite[d\'emonstration du th\'eor\`eme (4.3)]{WhG}). Il existe donc
une application $f\co M \to X$ induisant $\phi$. Par
transitivit\'e de l'action, il existe $h\in H^{1}(M;\Z/2\Z)$ tel
que
\begin{equation}
f^{*}(s_{X}) = h \cdot s_{M}. \label{eq:transitif}
\end{equation}
En appliquant $s \mapsto q_{s}$ \`a l'\'egalit\'e (\ref{eq:transitif}), nous
obtenons que
\begin{equation}
q_{f^{*}s_{X}} = q_{h \cdot s_{M}} = h \cdot q_{s_{M}}.
\end{equation}
Or
\begin{equation}
q_{f^{*}s_{X}} = f^{*}q_{s_{X}} = q_{s_{X}} \circ f^* = q_{s_{X}}
\circ \phi  = q_{s_{M}}.
\end{equation}
On en conclut que $q_{s_{M}} = h \cdot q_{s_{M}}$ d'o\`u $h = 1$. 
\endproof

\textbf{Remarque}\qua Le r\'esultat du Th\'eor\`eme
\ref{th:facteur-spin} reste vrai en rempla\c{c}ant spin par
spin$^c$. La d\'emonstration est essentiellement la m\^eme (en
utilisant \cite{DM}), la diff\'erence \'etant que les fonctions
quadratiques peuvent ne pas \^etre homog\`enes.
Voir \`a ce sujet \S \ref{sec:questions}, question 3.

Notons deux cons\'equences du Th\'eor\`eme \ref{th:facteur-spin}.
La premi\`ere utilise le r\'esultat principal de \cite{DM2}
relatif \`a la fonction de torsion de Turaev--Reidemeister
\cite{Turaev}. La fonction de Turaev--Reidemeister classifie les
structures spin$^c$ des lenticulaires \cite[\S 9.2]{Turaev0}. Le
r\'esultat suivant montre qu'elle est utile aussi dans l'\'etude
des applications de degr\'e un.

\begin{cor} \label{cor:torsion}
La fonction de torsion $T$ de Turaev-Reidemeister distingue
l'existence ou non d'une application de degr\'e un pr\'eservant
les structures spin (ou spin$^c$) d'une $3$-sph\`ere d'homologie
rationnelle $M$ sur un lenticulaire $L$.
\end{cor}

% attention a l'ordre des facteurs s_M s_X et la contra-variance...

% un petit exemple.

\`A l'aide du Th\'eor\`eme \ref{th:facteur-spin}, on montre \'egalement que le
Th\'eor\`eme  \ref{th:lenticulaire-proscrit} admet une version
parall\'elis\'ee. La v\'erification de ce fait est laiss\'ee au lecteur.

\section{D\'emonstrations des Th\'eor\`emes \ref{th:main} et
\ref{th:lenticulaire-proscrit}}

\label{sec:preuve-main}

\subsection{R\'esultats pr\'eliminaires} \label{subsec:prelim}

Notons ${\mathfrak{M}}$ le mono\"{\i}de des (classes
d'isomorphismes d') enlacements sur les 2-groupes. Pr\'esentons
tout d'abord quelques enlacements particuliers. Nous adoptons la
notation introduite dans \cite{KK}. Soit $k \geq 1$. Pour tout
entier impair $a$, on note $A^{k}(a)$ l'enlacement sur $\Z/2^k\Z$
qui envoie $(1\ {\rm{mod}}\ 2^k, \ 1\ {\rm{mod}}\ 2^k)$
 sur $\frac{a}{2^k}$ mod $1$. Sur $\Z/2^{k}\Z \times \Z/2^{k}\Z$, on d\'efinit
deux enlacements $E_{0}^{k}$ ($k \geq 1$) et $E_{1}^{k}$ ($k \geq 2$)
 comme suit:
$$\begin{array}{rccl}
E_{0}^{k}((x,y),(x',y')) & = & \frac{xy'+x'y}{2^{k}} &  {\rm{mod}}\ 1\\
E_{1}^{k}((x,y),(x',y')) & = & \frac{xx'+xy'+x'y+yy'}{2^{k}} &  {\rm{mod}}\ 1
\label{eq:E}
\end{array}$$
pour tous $x,y,x',y' \in \Z/2^{k}\Z$. Tout enlacement sur un
$2$-groupe fini est isomorphe \`a une somme ortho\-gonale
d'enlacements de type $A^{k}(a), E_{0}^{k}, $ et de $E_{1}^{k}$,
voir \cite{Wall}.

Le calcul suivant
\cite[Corollaire 2.2]{KK} est utile.

\begin{lem} \label{calcul}
Pour tout $k \geq 1$, $\sigma_{k}:{\mathfrak{M}} \to \overline{\Z}_{8}$ est
l'unique homomorphisme v\'erifiant les propri\'et\'es suivantes:
\begin{enumerate}
\item[$\bullet$] $\sigma_{k}(A^{l}(m)) = \left\{
\begin{array}{ll}
 (-1)^{\frac{m-1}{2}} & {\hbox{si}}\ l-k\ \hbox{is impair et positif}\\
m & {\hbox{if}}\ l-k \ {\hbox{is pair et positif}} \\
\infty  & {\hbox{si}}\ l=k  \\
0 & {\hbox{si}}\ l < k
\end{array} \right. $
\item[$\bullet$] $\sigma_{k}(E_{0}^{l}) = 0$, \item[$\bullet$]
$\sigma_{k}(E_{0}^{l}) = \left\{ \begin{array}{ll} 4 &
{\hbox{si}}\ l-k \ \hbox{is impair et positif}\\
0 & \hbox{sinon}. \end{array}\right.$
\end{enumerate}
\end{lem}

Nous avons donc, pour tout $k < l$, les \'egalit\'es suivantes dans $\Z/2\Z$:
$$\begin{array}{rclcl}
\sigma_{k}(A^{l}(m)) &  = & 1 \ {\rm{mod}}\ 2 & = & \rho_{k}(A^{l}(m)) \\
\sigma_{k}(E_{0}^{l})&  = &   0 \ {\rm{mod}}\ 2 & =&  \rho_{k}(E_{0}^l)\\
\sigma_{k}(E_{1}^{l})& = & 0 \ {\rm{mod}}\ 2&  = & \rho_{k}(E_{1}^l)
\end{array}$$
Comme les invariants $\sigma_{k}$ et $\rho_{k}$ sont additifs sur $\oplus$,
on en d\'eduit:

\begin{lem} \label{co2}
Soit $\lambda$ un enlacement sur un $\Z/2^l\Z$-module libre de type
fini. Pour tout $k < l$,
$\sigma_{k}(\lambda) = \rho_{l}(\lambda)\ {\rm{mod}}\ 2$.
\end{lem}

En vue du Th\'eor\`eme \ref{th:main}, nous avons besoin de la formule de
congruence modulo 4 suivante.

\begin{lem} \label{co3}
Soit $k, k+1$ deux \'el\'ements r\'eguliers pour un enlacement sur un
$\Z/2^l\Z$-module libre de type fini. Alors
\begin{equation}
\sigma_{k}(\lambda) + \sigma_{k+1}(\lambda) = 2\ \rho_{k}(\lambda)\
{\rm{mod}}\ 4. \label{eq:congruence4}
\end{equation}
\end{lem}

\proof[D\'emonstration]
D'apr\`es le lemme \ref{calcul}, $\sigma_{k}(A^{l}(m)) +
\sigma_{k}(A^{l}(m)) = m + (-1)^{\frac{m-1}{2}} \equiv 2 \
{\rm{mod}}\ 4$ (car $m$ est impair) et $\sigma_{k}(E_{0}^{l})
 + \sigma_{k}(E_{0}^{l}) \equiv \sigma_{k}(E_{1}^{l}) +
\sigma_{k}(E_{1}^{l}) \ {\rm{mod}}\ 4 \equiv 0 \ {\rm{mod}}\ 4$. 
\endproof

Remarquons une autre cons\'equence utile du lemme \ref{calcul}. Puisque les
facteurs orthogonaux cycliques sont les $A^{k}(m)$, le lemme \ref{calcul}
permet de d\'etecter la pr\'esence d'un tel facteur dans une d\'ecomposition
orthogonale:

\begin{lem} \label{detecte}
Un enlacement $\lambda$ admet un facteur orthogonal cyclique de
type $A^{k}(m)$ si et
seulement si $\sigma_{k}(\lambda) = \infty$.
\end{lem}

\subsection{D\'emonstration du Th\'eor\`eme \ref{th:main}}

\textbf{Necessit\'e}\qua Soit $\lambda$ un enlacement sur un $2$-groupe
$G$ fini. On v\'erifie que son tableau d'invariants est fini: soit $N$
l'exposant de $G$,
par d\'efinition, $\rho_{k}(\lambda) = 0$ d\`es que $k > N$; et nous avons
aussi $\sigma_{k}(\lambda) = 0$ d\`es que $k > N$ d'apr\`es le
lemme \ref{calcul}. Choisissons maintenant
une d\'ecomposition orthogonale $(G,\lambda) =
\oplus_{m} (G_{m}, \lambda_{m})$.
Soit $k$ un \'el\'ement r\'egulier pour le tableau d'invariants
associ\'e \`a $\lambda$.
Alors le facteur orthogonal
$(G_{k}, \lambda_{k})$ ne contient pas de facteur orthogonal cyclique
(c'est-\`a-dire un $A^{k}(m)$): sinon
$\sigma_{k}(\lambda_{k}) = \infty$ par le lemme \ref{detecte}, contredisant
le fait que $k$ est r\'egulier. Ainsi $(G_{k}, \lambda_{k})$ est une
somme orthogonale de copies de $E_{0}^{k}$ et de $E_{1}^{k}$. Puisque
$\rho_{k}(E_{0}^{k}) = \rho_{k}(E_{1}^{k}) = 2$, la condition $(1)$ en
r\'esulte. L'additivit\'e de $\sigma_{k}$ sur $\oplus$ et le lemme \ref{calcul}
impliquent $\sigma_{k}(\lambda) = \sum_{m \geq k+1} \sigma_{k}(\lambda_{m})$.
On en d\'eduit, avec le lemme \ref{co2}, la condition (2). Un argument tout
\`a fait similaire \`a partir du lemme \ref{co3} conduit \`a la condition (3).
V\'erifions \`a pr\'esent la condition (4). On observe tout d'abord que
tous les tableaux des types donn\'es ont une longueur impair. Aussi les
d\'elimiteurs $m < n$ v\'erifient $m = n\ \hbox{mod}\ 2$, de sorte que
 $\sigma_{m}(\lambda) - \sigma_{n}(\lambda) = \sum_{m<k<n} =
\sigma_{m}(\lambda_{k})$ d'apr\`es le lemme \ref{calcul}.  La suite
des v\'erifications pour la condition (4) est directe. 

Dans la d\'emonstration de la suffisance ci-dessous, si $C$ est un
entier positif ou nul et $\lambda$ un enlacement, afin d'all\'eger
les notations, on note
$C \cdot \lambda$ pour d\'esigner la somme orthogonale de $C$ copies
de l'enlacements $\lambda$.

\textbf{Suffisance}\qua La d\'emonstration se fait par
r\'ecurrence sur la longueur du tableau $T$. Soit $T$ un tableau
de longueur 1 v\'erifiant les conditions $(1)$ \`a $(4)$. Soit $m$
l'unique \'el\'ement de $I$. Si $m$ est r\'egulier, alors la
condition $(1)$ impose que $r(m) = 0$ mod $2$. On peut alors
prendre comme enlacement une somme orthogonale de $\frac{r(m)}{2}$
copies de $E_{0}^{m}$. Si $m$ n'est pas r\'egulier et $r(m) = 0$,
alors $s(m) = 0$ mod $2$ d'apr\`es la condition $(2)$ et donc on
peut prendre comme enlacement l'enlacement trivial. Si $s(m) =
\infty$, alors on peut prendre comme enlacement $A^{m}(1) \oplus
(r(m)-1)\ E_{0}^{k}$. Supposons avoir montr\'e qu'un tableau
$$T\co I= \{ m \in \N \ | \ m \geq k+1 \} \to \N \times
\overline{\Z}_{8}$$ satisfaisant les conditions $(1)$ \`a $(4)$
est admissible pour un enlacement $(G,\lambda) = \oplus_{l \geq k
+ 1} (G_{l}, \lambda_{l})$ o\`u chaque $(G_{l}, \lambda_{l})$ est
un enlacement sur un $\Z/2^l\Z$-module libre. Nous allons montrer
que tout tableau $T'\co \{k\} \cup I \to \N \times
\overline{\Z}_{8}$ qui prolonge $T$ et qui v\'erifie les
conditions $(1)$ \`a
$(4)$ est admissible.  

Si $k$ n'est pas r\'egulier alors on pose $\lambda_{k} = r(k) \cdot A^{k}(1)$.
On v\'erifie sans peine que $T'$ est admissible pour l'enlacement
$\lambda \oplus \lambda_{k}$. 

On suppose \`a pr\'esent que $k$ est r\'egulier. Par la condition (1),
$r(k) = 0$ mod $2$. Posons $\lambda_{k} = \frac{r(k)}{2} \cdot E_{0}^{k}$.
Clairement $\rho_{k}(\lambda_{k}) = r(k)$. Posons aussi $\lambda' =
\lambda_{k} \oplus \lambda$. Il y a trois cas \`a consid\'erer.

\textbf{Cas 1}\qua $k+1$ r\'egulier, $r(k+1) = 0$. Si
$k+2$ est r\'egulier, alors
\begin{align*}
s(k) & = s(k+2) && \text{par la condition (4)}\\
     & = \sigma_{k+2}(\lambda) && \text{d'apr\`es l'hypoth\`ese de r\'ecurrence sur $\lambda$}\\
     & = \sigma_{k}(\lambda) && \\
     & = \sigma_{k}(\lambda') && \text{d'apr\`es le Lemme \ref{calcul}.}
\end{align*} et $T'$ est admissible pour l'enlacement $\lambda'$. Si $k+2$
n'est pas r\'egulier, on applique la condition (3) au tableau $T'$
et \`a l'enlacement $\lambda'$ respectivement\footnote{On utilise
l'implication du Th\'eor\`eme \ref{th:main} d\'ej\`a d\'emontr\'ee
ci-dessus (``n\'ecessit\'e'') pour appliquer la condition (3) \`a
$\lambda'$.}:
\begin{equation}
\begin{array}{ccll}
s(k) + s(k+1) & \equiv & \displaystyle 2 \sum_{l \geq k+2} r(l) &
\hbox{mod}\ 4  \\
\sigma_{k}(\lambda') +  \sigma_{k+1}(\lambda') & \equiv &  \displaystyle 2
\sum_{l \geq k+2} \rho_{l}(\lambda') & \hbox{mod}\ 4.
\end{array} \label{eq:equa}
\end{equation}
Nous avons $\sigma_{k+1}(\lambda') = \sigma_{k+1}(\lambda_{k}) +
\sigma_{k+1}(\lambda) = \sigma_{k+1}(\lambda) = s(k+1)$.
Similairement, pour tout $l \geq k+1$, $\rho_{l}(\lambda') =
\rho_{l}(\lambda_{k}) + \rho_{l}(\lambda) = \rho_{l}(\lambda) =
r(l)$. Soustrayant l'une des \'egalit\'es \`a l'autre dans
$(\ref{eq:equa})$, on d\'eduit que $s(k) = \sigma_{k}(\lambda')\
{\hbox{mod}}\ 4$. Si $s(k) = \sigma_{k}(\lambda')$ mod $8$, alors
par d\'efinition, $T'$ est admissible pour $\lambda'$. Sinon $s(k)
= \sigma_{k}(\lambda') + 4$ mod $8$. Puisque
$\sigma_{k+2}(\lambda) = \infty$, par le lemme \ref{detecte},
$\lambda_{k+2}$ a un facteur ortho\-gonal cyclique qui est
$A^{k+2}(m)$ pour un certain entier impair $m$. Notons
$\lambda'_{k+2}$ le m\^eme enlacement que $\lambda_{k+2}$ mais en
rempla\c{c}ant ce facteur orthogonal cyclique par $A^{k+2}(m+4)$.
Il r\'esulte alors du lemme \ref{calcul} que $\sigma_{k+1}
(\lambda'_{k+2}) = \sigma_{k+1}(\lambda_{k+2})$. Par cons\'equent,
$T$ est aussi admissible pour l'enlacement $\lambda_{k+1} \oplus
\lambda'_{k+2} \oplus \displaystyle \bigoplus_{l \geq k+3}
\lambda_l$. Posons
$$\lambda'' = \lambda_{k} \oplus \lambda_{k+1} \oplus
\lambda_{k+2}' \oplus \displaystyle \bigoplus_{l \geq k+3}
\lambda_{l}.$$ On v\'erifie alors
\begin{align*}
s(k) & = \sigma_{k}(\lambda') + 4 && \\
     & = \sigma_{k}\Bigl(\lambda_k \oplus \lambda_{k+1} \oplus
     \lambda_{k+2} \oplus \bigoplus_{l \geq k+3}
     \lambda_{l}\Bigr) + 4 && \\
     & = \sigma_{k}(\lambda_{k}) + \sigma_{k}(\lambda_{k+1}) +
     \sigma_{k}(\lambda_{k+2})+ 4 + \sigma_{k}\Bigl(\bigoplus_{l \geq k+3}
     \lambda_{l}\Bigr) && \text{par additivit\'e de $\sigma_{k}$} \\
     & = \sigma_{k}(\lambda_{k}) + \sigma_{k}(\lambda_{k+1}) +
     \sigma_{k}(\lambda_{k+2}') + \sigma_{k}\Bigl(\bigoplus_{l \geq k+3}
     \lambda_{l}\Bigr) && \text{d'apr\`es le Lemme \ref{calcul}}\\
      & = \sigma_{k}(\lambda'') && \text{par additivit\'e de $\sigma_{k}$}
      \end{align*}
donc le tableau $T'$ est admissible pour $\lambda''$.

\textbf{Cas 2}\qua $k+1$ r\'egulier, $r(k+1) \not= 0$.
Le m\^eme argument que pr\'ec\'edemment donne $s(k) =
\sigma_{k}(\lambda')$ mod $4$. Si l'\'egalit\'e est vraie modulo
8, alors $T'$ est admissible pour $\lambda'$. Sinon $s(k) =
\sigma_{k}(\lambda') + 4$ et on proc\`ede de la fa\c{c}on
suivante. Puisque $\sigma_{k+1}(\lambda) \not= \infty$, par le
lemme \ref{detecte}, $\lambda_{k+1}$ n'a pas de facteur orthogonal
cyclique. Donc il existe $s, t \in \N$  tels que $\lambda_{k+1} =
s\ E_{0}^{k+1} \oplus t\ E_{1}^{k+1}$. D\'efinissons
$$\lambda_{k+1}' = \left\{ \begin{array}{ll} (s+1) \cdot E_{0}^{k+1}
\oplus (t-1)\cdot E_{1}^{k+1} & \hbox{si}\ t > 0 \\
(s-1) \cdot E_{0}^{k+1} \oplus E_{1}^{k+1} & \hbox{si}\ t = 0.
\end{array} \right.$$
Le Lemme \ref{calcul} implique que $\sigma_{k}(\lambda'_{k+1}) =
\sigma_{k}(\lambda_{k+1}) + 4$. De plus,
$\sigma_{k+1}(\lambda'_{k+1}) = \sigma_{k+1}(\lambda_{k+1}),$ de
sorte que $T$ est admissible pour l'enlacement $\lambda'_{k+1}
\oplus \displaystyle \bigoplus_{ l \geq k+2} \lambda_{l}$. Posons
$$\lambda'' = \lambda_{k} \oplus \lambda_{k+1}' \oplus\displaystyle
\bigoplus_{l \geq k+2} \lambda_{l}.$$ Alors une v\'erification
similaire \`a celle du cas pr\'ec\'edent montre que $s(k) =
\sigma_{k}(\lambda'')$. Ainsi le tableau $T'$ est admissible pour $\lambda''$.\\

\textbf{Cas 3}\qua $k+1$ n'est pas r\'egulier. Appliquons la
condition (2) \`a $T$ et $\lambda'$ respectivement:
\begin{equation}
 \begin{array}{ccll}
s(k) &  \equiv & \displaystyle \sum_{l \geq k+1} s(l) &
\hbox{mod}\ 2 \\
\sigma_{k}(\lambda') &  \equiv & \displaystyle \sum_{l \geq k+1}
\rho_{l}(\lambda') & \hbox{mod}\ 2 \\
\end{array} \label{eq:cond4applied}
\end{equation}
Nous avons $\sigma_{k+1}(\lambda') = \sigma_{k+1}(\lambda_{k}) +
\sigma_{k+1}(\lambda) = \sigma_{k+1}(\lambda) = s(k+1)$. De m\^eme, pour
$l \geq k+1$, $\rho_{l}(\lambda') = \rho_{l}(\lambda_{k}) + \rho_{l}
(\lambda) = \rho_{l}(\lambda) = r(l)$. Soustrayant l'une des \'egalit\'es
(\ref{eq:cond4applied}) \`a l'autre, on d\'eduit que $s(k) =
\sigma_{k}(\lambda')$ mod $2$. Si cette derni\`ere \'egalit\'e reste
vraie modulo $8$, alors $T'$ est admissible pour $\lambda'$. Si $s(k)
= \sigma_{k}(\lambda') \pm 2$, on proc\`ede de la fa\c{c}on suivante.
Puisque $\sigma_{k+1}(\lambda) = \infty$, d'apr\`es le lemme \ref{detecte},
$\lambda_{k+1}$ a un facteur orthogonal cyclique $A^{k+1}(m)$,
o\`u $m$ est un entier impair. D\'efinissons $\lambda'_{k+1}$ comme \'etant le m\^eme
enlacement que $\lambda_{k+1}$ mais en rempla\c{c}ant ce facteur par $A^{k+1}(m+s(k)-
\sigma_{k}(\lambda')) = A^{k+1}(m\pm 2)$. Alors $\sigma_{k}(\lambda'_{k+1}) =
\sigma_{k}(\lambda_{k+1})$. Posons $\lambda'' =\lambda_{k} \oplus \lambda'_{k+1} \oplus
\displaystyle \bigoplus \lambda_{l}$. Alors $s(k) = \sigma_{k}(\lambda'')$ et le
tableau $T'$ est admissible pour $\lambda''$. Ainsi il ne reste \`a consid\'erer que le cas
o\`u $s(k) = \sigma_{k}(\lambda') + 4$. Il y a deux possibilit\'es.

\textbf{Possibilit\'e 1}\qua $s(k+2) = \infty$. Le lemme \ref{detecte} dit alors que
$\lambda_{k+2}$ a un facteur orthogonal cyclique $A^{k+2}(m)$ pour un certain entier impair
$m$. Rempla\c{c}ons le par $A^{k+2}(m+4)$ et renommons le nouvel enlacement $\lambda_{k+2}'$.
Une v\'erification similaire \`a celle du Cas $1$ montre que $T'$ est admissible pour
l'enlacement $\lambda'' = \lambda_{k} \oplus \lambda_{k+1} \oplus \lambda_{k+2}'
\oplus \displaystyle \bigoplus_{l \geq k+3} \lambda_{l}$.

\textbf{Possibilit\'e 2}\qua $s(k+2) \not= \infty$. Alors $r(k+1) \geq 2$: en effet,
sinon $r(k+1) = 1$ ($k+1$ n'est pas r\'egulier) et nous appliquons la condition $(4)$
\`a $T$ et $\lambda'$ respectivement
\begin{equation}
 \begin{array}{ccl}
s(k) - s(k+2) & = &\pm 1 \\
\sigma_{k}(\lambda') - \sigma_{k+2}(\lambda') & =& \pm 1
\end{array} \label{eq:c6}
\end{equation}
et on d\'eduit $s(k) - \sigma_{k}(\lambda') = 0$ ou $\pm 2$, une contradiction. On traite alors
les deux cas s\'epar\'ement.

\begin{enumerate}
\item[-] Si $r(k+1) \geq 3$, alors on affirme que $\lambda_{k+1}$ a un facteur orthogonal
$S = E_{0}^{k+1}$ ou $E_{1}^{k+1}$. [Preuve: sinon $\lambda_{k+1}$ a au moins
trois facteurs orthogonaux cycliques $A^{k+1}(n_1)$, $A^{k+1}(n_2)$ et
$A^{k+1}(n_3)$. Les relations $(0.2)$ et $(0.3)$ de \cite{KK} impliquent
alors que les $n_i$ sont deux \`a deux distincts dans $\{ \pm 1, \pm 5 \} =
(\Z/8\Z)^{\times} = \Z/2\Z \times \Z/2\Z$.
Il s'ensuit qu'il existe $i, j$ tels que $n_j = n_i + 4$ mod $8$. Apr\`es
renum\'erotation, on peut supposer $i = 1, j = 2$. Mais alors, d'apr\`es
\cite[rel. (0.1)]{KK} $A^{k+1}(n_1) \oplus A^{k+1}(n_3) = A^{k+1}(n_2) \oplus A^{k+1}(n_3 + 4)$ et donc
$A^{k+1}(n_1) \oplus A^{k+1}(n_2) \oplus A^{k+1}(n_3) =
2 \ A^{k+1}(n_2) \oplus A^{k+1}(n_3)$ et les deux premiers indices
sont \'egaux \`a $n_2$,
contradiction.] Nous posons alors
$$ S' = \left\{ \begin{array}{ll} E_{1}^{k+1} & \hbox{si}\ S = E_{0}^{k+1}\\
E_{0}^{k+1} & \hbox{si}\ S = E_{1}^{k+1} \end{array}\right.$$
D\'esignons par $\lambda_{k+1}'$ l'enlacement $\lambda_{k+1}$ o\`u
l'on a remplac\'e $S$ par $S'$. Alors $\sigma_{k+1}(\lambda_{k+1}) = \sigma_{k+1}(\lambda_{k+1}) + 4$.
On conclut, comme dans le Cas 2, que $T'$ est admissible pour l'enlacement $\lambda'' = \lambda_{k} \oplus
\lambda_{k+1}' \oplus \displaystyle \bigoplus_{l \geq k+2} \lambda_{l}$.
\item[-] Si $r(k+1) = 2$, vu que d'apr\`es le lemme \ref{detecte}, $\lambda_{k+1}$ admet d\'ej\`a un
facteur orthogonal cyclique, la seule possibilit\'e est
$\lambda_{k+1} = A^{k+1}(m) \oplus A^{k+1}(n)$ pour des \'el\'ements inversibles $m, n \in \Z/8\Z$. Nous
affirmons que $m = n$ mod $4$. [Sinon, le lemme \ref{calcul} donne $\sigma_{k}(\lambda') -
\sigma_{k+2}(\lambda') = 0$ alors que la condition $(4)$ appliqu\'ee \`a $T'$ implique $s(k) - s(k+2)
 = 0$ ou $\pm 2$. Soustrayant la premi\`ere \'egalit\'e \`a la seconde, on trouve
$s(k) - \sigma_{k}(\lambda') = 0$ ou $\pm 2$, contradiction.] Rempla\c{c}ons alors dans $\lambda_{k+1}$
les facteurs $A^{k+1}(m)$ et $A^{k+1}(n)$ par $A^{k+1}(m+4)$ et $A^{k+1}(n)$ respectivement. Notons
$\lambda_{k+1}'$ le nouvel enlacement qui en r\'esulte. Alors d'apr\`es le lemme \ref{calcul},
$\sigma_{k}(\lambda_{k+1}) - \sigma_{k}(\lambda_{k+1}') = 2 \ ( (-1)^{\frac{m-1}{2}} +
(-1)^{\frac{n-1}{2}})  = 4$ mod $8$. Il s'ensuit que $\sigma_{k}(\lambda_{k+1}') - s(k) = 0$.
D\'efinissons $\lambda'' = \lambda_{k} \oplus \lambda_{k+1}' \oplus \displaystyle \bigoplus_{l \geq k+2}
 \lambda_{l}$. Alors $s(k) = \sigma_{k}(\lambda'')$  et $T'$ est admissible pour $\lambda''$.\endproof

\end{enumerate}

\textbf{Remarque}\qua La d\'emonstration est constructive. Si l'on
suppose construit
\`a l'\'etape $k$ l'enlacement, alors on peut le construire \`a l'\'etape $k+1$. Il est possible
de raffiner la construction de sorte que l'enlacement obtenu soit sous la forme normale
d\'ecrite dans \cite[\S\S 3, 4]{Mi}.

% dire poruquoi c'est constructif. Si l'on suppose construit
% a l'etape k l'enlacement, alors on peut le construire a l'etape k+1.

\subsection{D\'emonstration du Th\'eor\`eme \ref{th:lenticulaire-proscrit}}
\label{subsec:demth6}

On note $(\Z/8\Z)^{\times} = \{\pm 1, \pm 3\ {\rm{mod}}\ 8\}$ le
groupe des \'el\'ements inversibles de l'anneau $\Z/8\Z$. Il est
isomorphe \`a $\Z/2\Z \times \Z/2\Z$. Il est bien connu que la
classe d'isomorphisme de l'enlacement d'un lenticulaire $L(q,q')$
ne d\'epend que du r\'esidu quadratique de $q'$ modulo $q$. En
particulier, si $q$ est une puissance de $2$, elle ne d\'epend que
$q$ modulo $8$ (voir par exemple \cite[p.69]{Vino}). L'application
$(\Z/8\Z)^{\times} \times (\Z/8\Z)^{\times} \times
(\Z/8\Z)^{\times} \to (\Z/8\Z)^{\times}$ d\'efinie par $$ (r_{1},
r_{2}, r_{3}) \mapsto s = 4 - r_{2} + (-1)^{\frac{r_{1}+1}{2}} +
(-1)^{\frac{r_{3}+1}{2}}\ {\rm{mod}}\ 8$$ est surjective.
Consid\'erons la $3$-vari\'et\'e $M = L(8,r_{1}) \# L(16,r_{2}) \#
L(32, r_{3})$. Nous affirmons que $M$ n'admet pas d'application de
degr\'e $1$ sur $L(16,s)$ alors qu'elle admet une application de
degr\'e $1$ sur chaque lenticulaire $L(16,r)$ pour $r \not\equiv
s$ mod $8$. Pour le voir, posons $\lambda = \lambda_{M}$ et
$\lambda' = A^{5}(r)$ (l'enlacement cyclique sur $\Z/32\Z$
envoyant $(1,1)$ sur $r/32$ mod $1$). D'apr\`es le Th\'eor\`eme
\ref{th:facteur}, nous avons \`a voir pour quel $r \in
(\Z/8\Z)^{\times}$ il existe un tableau admissible dans
$S_{\lambda, \lambda'}$. Puisque $\sigma_{m}(\lambda') = \infty$
si et seulement si $m = 4$, il y a $9$ tableaux distincts dans
$S_{\lambda, \lambda'}$. Tout tableau $T \in S_{\lambda,
\lambda'}$ contient un sous-tableau extrait de la forme
$$\begin{array}{|c|c|c|c|c|} \hline
1 & 2 & 3 & 4 & 5 \\
\hline \hline
0 & 0 & 1 & 0 & 1 \\
\hline
\sigma_{1}(\lambda) - \sigma_{1}(\lambda') & \sigma_{2}(\lambda) -
\sigma_{2}(\lambda') & \infty & x & \infty \\
\hline
\end{array}$$ o\`u $x \in \overline{\Z}_{8}$. (Le tableau $T$ entier s'obtient en
prolongeant trivialement le tableau ci-dessus.) D'apr\`es les
calculs de \S\S \ref{subsec:prelim}, nous avons
\begin{align*}
\sigma_{1}(\lambda) - \sigma_1(\lambda') & =
r_{1}+(-1)^{\frac{r_{2}-1}{2}}+r_{3}-(-1)^{\frac{r-1}{2}}\
{\rm{mod}}\ 8 \\
{\rm{et}}\ \ \sigma_{2}(\lambda) - \sigma_{2}(\lambda') & =
 (-1)^{\frac{r_1 - 1}{2}} + r_2 + (-1)^{\frac{r_3 - 1}{2}} - r\ {\rm{mod}}\
 8.
\end{align*}
Examinons alors les conditions du Th\'eor\`eme \ref{th:main}: il
est ais\'e de constater qu'elles sont toutes remplies si et
seulement si la derni\`ere condition $(4)$ est remplie,
c'est-\`a-dire si et seulement si $$x \equiv \pm 1\ {\rm{mod}}\ 8\
\ \ {\rm{et}}\ \ (\sigma_{2}(\lambda) - \sigma_{2}(\lambda')) - x
\equiv \pm 1\ {\rm{mod}}\ 8.$$ (Ces deux conditions correspondent
aux deux sous-tableaux distingu\'es de type T$_{1}$ que l'on peut
extraire du tableau $T$. Il s'agit des sous-tableaux
$\begin{array}{|c|} \hline 3
\\ \hline \hline 1 \\ \hline \infty
\\ \hline
 \end{array}$ et  $\begin{array}{|c|} \hline 5 \\ \hline \hline
1 \\ \hline \infty \\ \hline
 \end{array}$ respectivement.)
  Par
cons\'equent, $T$ est admissible si et seulement s'il existe $r
\in (\Z/8\Z)^{\times}$ tel que $$(-1)^{\frac{r_1 - 1}{2}} + r_2 +
(-1)^{\frac{r_{3}-1}{2}} - r = \pm 1 \pm 1  \in \{ 0, \pm 2 \}
\subset \Z/8\Z.$$
Le seul cas o\`u cette condition est mise en d\'efaut (pour $r \in
(\Z/8\Z)^{\times}$) est quand $r = s(r_1, r_2, r_3) = 4 - r_{2} +
(-1)^{\frac{r_{1}+1}{2}} + (-1)^{\frac{r_{3}+1}{2}}$ mod $8$ comme
d\'efini plus haut. Notre affirmation en r\'esulte. Pour finir,
d'apr\`es le lemme \ref{lem:chirurgie}, on peut chirurgiser $M$ de
mani\`ere \`a rendre $M$ irr\'eductible (en fait hyperbolique) et
\`a pr\'eserver les
propri\'et\'es ci-dessus. \endproof

\begin{lem} \label{lem:chirurgie}
Soit $X$ une vari\'et\'e connexe ferm\'ee orient\'ee de dimension
$3$. Il existe une infinit\'e de vari\'et\'es $M$ irr\'eductibles
(hyperboliques) ayant la m\^eme alg\`ebre de cohomologie et m\^eme
enlacement que $X$ et admettant une application $M \to X$ de
degr\'e $1$.
\end{lem}

\proof[D\'emonstration] D'apr\`es un r\'esultat de Myers
\cite{My}, $X$ contient un n\oe{u}d nul-homotope $K$. Ceci
implique que le compl\'ement d'un voisinage r\'egulier de $K$ est
irr\'eductible (hyperbolique). Le c\'el\`ebre th\'eor\`eme de
Gordon--Luecke \cite{GL} dit alors qu'il existe un nombre infini
de remplissages (in\'equivalents) de Dehn sur $K$ produisant
chacun une vari\'et\'e hyperbolique $M$ avec la m\^eme alg\`ebre
de cohomologie et enlacements isomorphes. Un argument d\^u \`a
Boileau--Wang \cite[preuve de la Prop. 3.2]{BW} construit
explicitement une application $M \to X$
de degr\'e $1$. \endproof

\section{Quelques questions} \label{sec:questions}

Nous incluons dans cette section quelques questions sugg\'er\'ees par
les r\'esultats de cet article.

\noindent{\bf{Question 1}}\qua [Structure du mono\"{\i}de des enlacements]\qua
Calculer le nombre de classes d'isomorphismes d'enlacements de
rangs donn\'es. Plus g\'en\'eralement, avec les notations introduites
\`a la fin de  \S \ref{sec:monoide}, d\'ecrire les fibres de l'application
Profil. Plus pr\'ecis\'ement, existe-t-il une classe ${\mathcal{C}}$
de parties de $\N$ telles que ${\mathfrak{M}}$ soit classifi\'e
\`a partir des fibres Profil$^{-1}(E)$, $E \in {\mathcal{C}}$ ? On peut
poser les m\^emes questions en un sens plus faible en rempla\c{c}ant les
fibres par leur groupe de Witt. Enfin, ces
questions peuvent \^etre pos\'ees au sujet du mono\"{\i}de des
 formes quadratiques.

\noindent{\bf{Question 2}}\qua [R\'ealisabilit\'e des enlacements par
les vari\'et\'es de Seifert]\qua Il est connu \cite[Th.~6.1]{KK} que
tout enlacement peut \^etre r\'ealis\'e comme l'enlacement d'une
vari\'et\'e connexe ferm\'ee orient\'ee de dimension $3$. On peut
modifier l'argu\-ment du lemme \ref{lem:chirurgie} pour imposer
que la $3$-vari\'et\'e r\'ealisant l'enlacement soit
irr\'eductible. Est-ce que l'on peut imposer que la
$3$-vari\'et\'e soit de Seifert (une question qui a peut-\^etre
motiv\'e Seifert \`a introduire les vari\'et\'es de Seifert) ?
Peut-\^etre peut-on
le d\'emontrer \`a l'aide des techniques de \cite{BW}. 

\noindent{\bf{Question 3}}\qua [Forme normale]\qua L'article \cite{Mi} contient
un algorithme de mise sous forme normale pour tout enlacement et permet
de simplifier la preuve originale de \cite{KK} de la pr\'esentation par
g\'en\'erateurs et relations de ${\mathfrak{M}}$. Peut-on g\'en\'eraliser
le r\'esultat principal (Th.~4.4) de \cite{Mi} aux formes quadratiques sur
les $2$-groupes ?

\noindent{\bf{Question 4}}\qua [Reconnaissance de fonctions quadratiques]\qua
 D'apr\`es
la remarque 2 \`a la fin de \S \ref{sec:reconnaissance}, le Th\'eor\`eme
\ref{th:classification-formesquad} est une g\'en\'eralisation du Th\'eor\`eme
\ref{th:main} aux
{\emph{formes}} quadratiques. Quoique plus d\'elicat, il serait
utile de les g\'en\'eraliser aux
{\emph{fonctions}} quadratiques, c'est-\`a-dire aux applications
$q:G \to \Q/\Z$ telle que $b_{q}(x,y) = q(x+y) - q(x) - q(y)$ soit
bilin\'eaire en $x$ et $y$ {\emph{sans requ\'erir a priori}} la condition
d'homog\'en\'eit\'e $q(n\ x) = n^2\ q(x)$. Une telle question est motiv\'ee,
par exemple, par le fait qu'une $3$-vari\'et\'e \'equip\'ee d'une structure
Spin$^c$ poss\`ede canoniquement une telle fonction quadratique. (C'est
d'ailleurs ce qui permet de g\'en\'eraliser le Th\'eor\`eme
\ref{th:facteur-spin} \`a ce cadre.) Ce fait
est vrai plus g\'en\'eralement d'ailleurs, pour toute vari\'et\'e
ferm\'ee de dimension $4n-1$ \'equipp\'ee d'une structure complexe sur
$TM \oplus \R_M$ (fibr\'e tangent stabilis\'e une fois) \cite{LW}.

\Addresses\recd

\end{document}